\documentclass[10pt]{amsart}
\usepackage{amssymb,latexsym,amsmath}
\usepackage{graphics}                 
\usepackage{color}                    
\usepackage{hyperref}                 
\usepackage[all]{xy}

\begin{document}

\newtheorem{definition}{Definition}
\newtheorem{theorem}{Theorem}[section]
\newtheorem{proposition}[theorem]{Proposition}
\newtheorem{lemma}[theorem]{Lemma}
\newtheorem{corollary}[theorem]{Corollary}
\newtheorem{question}[theorem]{Question}
\newtheorem{remark}[theorem]{Remark}
\newtheorem{example}[theorem]{Example}
\newtheorem{conjecture}[theorem]{Conjecture}

\newtheorem{correction}{Correction}

\def\A{{\mathbb{A}}}
\def\C{{\mathbb{C}}}
\def\H{{\mathbb{H}}}
\def\L{{\mathbb{L}}}
\def\P{{\mathbb{P}}}
\def\Q{{\mathbb{Q}}}
\def\Z{{\mathbb{Z}}}
\def\Ch{{\rm Ch}}
\def\p{{\mathbf{p}}}

\def\id{{\rm id}}
\def\sp{{\rm SP}}
\def\mC{{\mathcal{C}}}

\title{The Lawson-Yau Formula and its generalization}
\author{Wenchuan Hu }

\address{
Department of Mathematics, MIT, Room 2-363B,
77 Massachusetts Avenue, Cambridge, MA 02139, USA }
\email{wenchuan@math.mit.edu}
\keywords{Algebraic cycles; Chow varieties; Euler characteristic, Group action}
\date{December 5, 2008}
\maketitle
\pagestyle{myheadings}
 \markright{The Lawson-Yau formula and its generalization}

\begin{abstract}
The Euler characteristic of Chow varieties of algebraic cycles of a given degree
in complex projective spaces was computed by Blaine Lawson and Stephen Yau by
using holomorphic symmetries of cycles spaces.  In this paper we compute this
in a direct and elementary way and generalize this formula to the $l$-adic
 Euler-Poincar\'{e} characteristic for Chow varieties over any algebraically closed field.
 Moreover, the Euler characteristic for Chow varieties with certain  group action is
 calculated. In particular, we calculate the Euler characteristic of the space of
  \emph{right} quaternionic cycles of a given
 dimension and degree in complex projective spaces.
\end{abstract}
\tableofcontents

\section{Introduction}

Let $\P^n$ be the  complex projective space of  dimension $n$ and let $C_{p,d}(\P^n)$
be the space of effective algebraic $p$-cycles of degree $d$ on $\P^n$. A fact proved
by Chow and Van der Waerden  is that $C_{p,d}(\P^n)$ carries the structure of a
closed complex algebraic variety \cite{Chow-Waerden}.
Hence it carries
the structure of a compact Hausdorff space.
Denote by $\chi(C_{p,d}(\P^n))$ the Euler Characteristic of $C_{p,d}(\P^n)$.
 This number was computed in terms of $p, d$ and $n$ by Blaine Lawson and Stephen Yau explicitly, i.e.,

\begin{theorem}[Lawson-Yau, \cite{Lawson-Yau}]
\begin{equation}\label{eq:LY}
\chi(C_{p,d}(\P^n))=\big(^{v_{p,n}+d-1}_{\quad\quad d}\big), \quad\hbox{where $v_{p,n}=(^{n+1}_{p+1})$}.
\end{equation}
\end{theorem}

This equation (\ref{eq:LY}) is called the \textbf{Lawson-Yau formula}.
 The original method of calculation  is an application of a fixed point formula
 for a compact complex analytic space with a weakly holomorphic $S^1$-action.

Equivalently, if we define
$$ Q_{p,n}(t):=\sum_{d=0}^{\infty} \chi(C_{p,d}(\P^n)) t^d,
$$
then the Lawson-Yau formula may be restated as
\begin{equation}\label{eq1.1}
Q_{p,n}(t)=\bigg(\frac{1}{1-t}\bigg)^{(^{n+1}_{p+1})}, \quad\hbox{where $\chi(C_{p,0}(\P^n)):=1$} .
\end{equation}

A direct and elementary  proof of  the Lawson-Yau formula is given in section \ref{sec2}
by using the technique ``pulling of normal cone" established in the book of
Fulton (\cite{Fulton}), which was used by Lawson in proving his Complex
Suspension Theorem (\cite{Lawson1}).

As an application of this elementary method, we obtain an $l$-adic version of the Lawson-Yau formula:
\begin{theorem} \label{main}
Let $C_{p,d}(\P^n)_K$ be the space of effective $p$-cycles of degree $d$ in $\P^n_K$.  For all $l$ prime to $char(K)$,
we have
\begin{equation}\label{eqn1.3}
\chi(C_{p,d}(\P^n)_K,l)=\big(^{v_{p,n}+d-1}_{\quad\quad d}\big), \quad\hbox{where $v_{p,n}=(^{n+1}_{p+1})$},
\end{equation}
where $\chi(X_K,l)$ denotes the $l$-adic Euler-Poincar\'{e} Characteristic of an algebraic variety $X_K$ over $K$.
\end{theorem}

The detailed explanation of notations in Theorem \ref{main} as well as its proof is given in section \ref{sec3}.

We apply our method to the space of algebraic cycles with certain finite group action $G$
to obtain the Euler characteristic of the $G$-invariant Chow varieties (see Theorem \ref{Th4.4}).

As an application, we obtain the Euler characteristic of the space of \emph{right} quaternionic cycles of a given
 dimension and degree in complex projective spaces (see Corollary \ref{cor5.2} ).

\thanks{\emph{Acknowledgements}: I would like to gratitude Michael Artin,
 Eric Friedlander and James McKernan for their interest and encouragement
as well as helpful advice on the organization of the paper.
}

\section{An elementary proof of the Lawson-Yau formula}\label{sec2}

First we list a result we will use in our calculation.  See \cite{Fulton2}, page 95.
\begin{lemma}\label{lemma2.1}
Let $X$ be a complex algebraic variety (not necessarily smooth, compact, or irreducible)
 and let $Y\subset X$ be a closed algebraic set with complement $U$. Then $\chi(X)=\chi(U)+\chi(Y)$.
\end{lemma}

Now we give a brief review of Lawson's construction of the suspension map and some of its properties.
The reader is referred to \cite{Lawson1} for details.
Fix a hyperplane $\P^n\subset \P^{n+1}$ and a point $P=[0:\cdots:0:1]\in \P^{n+1}-\P^n$.
 Let $V\subset \P^n$ be any closed algebraic subset. The \textbf{algebraic suspension
of $V$} with vertex $P$ (i.e., cone over $P$) is the set
$$\Sigma_P V:=\bigcup\{ l~|~ l \hbox{ is a projective line through $P$ and intersects $V$}\}.
$$
Extending by linearity  algebraic suspension gives a continuous homomorphism
\begin{equation}\label{eqn2}
 \Sigma_P: C_{p,d}(\P^n)\to C_{p+1,d}(\P^{n+1})
\end{equation}
for any $p\geq 0$.

Set $$T_{p+1,d}(\P^{n+1}):=\bigg\{ c=\sum n_iV_i\in C_{p+1,d}(\P^{n+1})| \dim(V_i\cap \P^n)=p, \forall i\bigg\}$$
and $B_{p+1,d}(\P^{n+1}):=C_{p+1,d}(\P^{n+1})-T_{p+1,d}(\P^{n+1})$.

The following proposition proved by Lawson in \cite{Lawson1} is the key point in our calculation.
An algebraic version  was given by Friedlander in \cite{Friedlander1}.
\begin{proposition}\label{deform}
 The subset $T_{p+1,d}(\P^{n+1})\subset T_{p+1,d}(\P^{n+1})$ is Zariski open.
 Moreover, the image $\Sigma_P: C_{p,d}(\P^n)\to C_{p+1,d}(\P^{n+1})$ is included in $T_{p+1,d}(\P^{n+1})$
 and $ \Sigma_P (C_{p,d}(\P^n))\subset T_{p+1,d}(\P^{n+1})$ is a strong deformation retract.
\end{proposition}

In particular, their Euler characteristics coincide, i.e., we have
\begin{equation}\label{eqn3}
 \chi(C_{p,d}(\P^n))=\chi(T_{p+1,d}(\P^{n+1})).
\end{equation}

By Proposition \ref{deform}, $B_{p+1,d}(\P^{n+1})$ is a closed subset of $C_{p+1,d}(\P^{n+1})$. From the definition,
$$B_{p+1,d}(\P^{n+1})=\{ c=\sum n_iV_i\in C_{p+1,d}(\P^{n+1})| V_i\subset \P^n, \hbox{for some $i$}\},$$ i.e.,
 there is at least one irreducible component lying in the  fixed hyperplane $\P^n$.

\begin{lemma}\label{lemma2.2}
 $B_{p+1,d}(\P^{n+1})=\amalg_{i=1}^{d}B_{p+1,d}(\P^{n+1})_i$ , where $\amalg$ means disjoint union and
 $$B_{p+1,d}(\P^{n+1})_i=\left\{
                 c\in B_{p+1,d}(\P^{n+1})
		 \begin{array}{ll}\left|\begin{array}{ll}
		  &c=\sum n_kV_k+\sum m_jW_j,   \\
                  &  V_k\subset \P^n,\forall k,\\
		   & \dim(W_j\cap \P^n)=p, \forall j\\
		   & \deg(\sum n_kV_k)=i,\\
		   & and ~~\deg(\sum m_jW_j)=d-i.
		   \end{array}
                  \right.
                \end{array}
              \right\}.
$$


 For each $i$, $B_{p+1,d}(\P^{n+1})_i=C_{p+1, i}(\P^n)\times T_{p+1,d-i}(\P^{n+1})$.
\end{lemma}
\begin{proof} Clear from the definition of $B_{p+1,d}(\P^{n+1})_i$.
\end{proof}

From Lemma \ref{lemma2.2}, we have $$\chi(B_{p+1,d}(\P^{n+1})_i)=\chi(C_{p+1, i}(\P^n))\cdot\chi(T_{p+1,d-i}(\P^{n+1})).$$ Hence we get
$$
 \begin{array}{llll}
\chi(B_{p+1,d}(\P^{n+1}))
                 &=&\sum_{i=1}^d \chi(B_{p+1,d}(\P^{n+1})_i)\\
                 &&\hbox{(by inclusion-exclusion principle)}\\
                 &=&\sum_{i=1}^d \chi(C_{p+1, i}(\P^n))\cdot\chi(T_{p+1,d-i}(\P^{n+1}))\\
                 &=& \sum_{i=1}^d \chi(C_{p+1, i}(\P^n))\cdot\chi(C_{p,d-i}(\P^{n})), \\ &&\hbox{(by equation (\ref{eqn3}))}\\
                \end{array}
$$

Therefore we have the following result:

\begin{proposition}\label{prop2.2}
 For any integer $p\geq 0$ and $d\geq 1$, we have the following recursive formula
\begin{equation}\label{eqn4}
\chi(C_{p+1, d}(\P^{n+1}))=\chi(C_{p, d}(\P^{n}))+\sum_{i=1}^d \chi(C_{p+1, i}(\P^n))\cdot\chi(C_{p,d-i}(\P^{n})),
\end{equation}
where $\chi(C_{q, 0}(\P^{N}))=1$ for integers $N\geq  q\geq 0$. In particular,
when $d=1$,  equation (\ref{eqn4}) is just the combinatorial identity $(^{n+2}_{p+2})=(^{n+1}_{p+1})+(^{n+1}_{p+2})$.
\end{proposition}

To compute $\chi(C_{p, d}(\P^{n}))$, it is enough to identify the initial values.

\begin{lemma}\label{lemma2.3}
$\chi(C_{0, d}(\P^{n}))=(^{n+d}_{~~d})$.
\end{lemma}
 The equality is a special case of MacDonald formula (\cite{Macdonald}).

\begin{proof}[Proof of Lemma \ref{lemma2.3}]
Now we give an independent proof for MacDonald formula in this special case.
  We can write $$C_{0,d}(\P^{n+1})=C_{0,d}(\C^{n+1})\coprod B_{0,d}(\P^{n+1}),$$
  where $C_{0,d}(\C^{n+1})\subset C_{0,d}(\P^{n+1})$ contains effective 0-cycles $c$
  of degree $d$ such that no points in $c$ lying in the fixed hyperplane $\P^n$ and
   $B_{0,d}(\P^{n+1})$ is the complement of $C_{0,d}(\C^{n+1})$ in  $C_{0,d}(\P^{n+1})$.
   It is easy to see that $C_{0,d}(\C^{n+1})$ is contractible.
   We can write $B_{0,d}(\P^{n+1})=\coprod_{i=1}^d B_{0,d}(\P^{n+1})_i$ as in Lemma \ref{lemma2.2}, where
 $B_{0,d}(\P^{n+1})_i$ contains 0-cycles $c$ of degree $d$
 on $\P^{n+1}$ in which there are exact $i$ points (count multiplicities) lying in $\P^n$, hence
 $B_{0,d}(\P^{n+1})_i=C_{0,i}(\P^n)\times C_{0,d-i}(\C^{n+1})$. In particular,
 $\chi(B_{0,d}(\P^{n+1}))=\sum_{i=1}^d\chi(C_{0,i}(\P^n))$.

 Therefore, we have
 $$\chi(C_{0, d}(\P^{n+1}))=1+\sum_{i=1}^d\chi(C_{0,i}(\P^n)).$$
The first formula in the lemma follows from this by induction.
\end{proof}

\begin{proof}[Proof of Theorem \ref{main}]
Equation (\ref{eqn4}) together with Lemma \ref{lemma2.3} is equivalent to the following recursive functional equation with initial values
\begin{equation}\label{eqn5}
\begin{array}{lll}
Q_{p+1,n+1}(t)&=&Q_{p+1,n}(t)\cdot Q_{p,n}(t),\\
Q_{0,m}(t)&=&(\frac{1}{1-t})^{m+1}.
\end{array}
\end{equation}

From this, we get the equation (\ref{eq1.1})  by induction on $n$ and hence the Lawson-Yau formula (\ref{eq:LY}).
\end{proof}

\begin{example}
For divisors of degree $d$ in $\P^n$, we have the formula $\chi(C_{p, d}(\P^{p+1}))=(^{p+d+1}_{\quad d})$.
\end{example}
 From equation \ref{eqn5}, we have
 $$
 \begin{array}{lll}
 Q_{p,p+1}(t)&=&Q_{p,p}(t)\cdot Q_{p-1,p}(t)\\
 &=& \frac{1}{1-t}\cdot Q_{p-1,p}(t)\\
  \end{array}
 $$
since  $C_{p,d}(\P^p)$ contains exactly one degree $d$ cycle and so $\chi(C_{p,d}(\P^p))=1$, i.e, $Q_{p,p}(t)= \frac{1}{1-t}$.
By the fact that $Q_{0,1}(t)=(\frac{1}{1-t})^{2}$ and induction on $p$, we get  $ Q_{p,p+1}(t)=(\frac{1}{1-t})^{p+2}$.
Hence $\chi(C_{p, d}(\P^{p+1}))=(^{p+d+1}_{\quad d})$.

 Alternatively, this formula follows directly from the fact that $C_{p, d}(\P^{p+1})$ is
 the moduli space of hypersurfaces of degree $d$ in $\P^{p+1}$ and hence it
is a complex projective space of dimension $(^{p+d+1}_{\quad d})-1$. To see this,  we choose
a basis for the monomials of degree $d$ in $p+2$ variables and then associate a point in this
projective space to the hypersurface whose defining equation is  given by the coordinates of that point (cf. \cite{Friedlander1}).

\section{The $l$-adic Euler-Poincar\'{e} Characteristic for Chow varieties}\label{sec3}.
In the section, the Lawson-Yau formula is generalized to an algebraically closed field
 $K$ with arbitrary characteristic $char(K)\geq 0$. Let $l$ be a positive integer prime to $char(K)$.

For a variety $X$ over $K$, let $H^i(X,\Z_l)$ be the $l$-adic cohomology group of $X$.
 Set $H^i(X,\Q_l):=H^i(X,\Z_l)\otimes_{\Z_l} \Q_l$. Denote by $\beta^i(X,l):=\dim_{\Q_l}H^i(X,\Q_l)$
 the  $i$-th $l$-adic Betti number of $X$. The $l$-adic Euler Characteristic is defined by
 $\chi(X,l):=\sum_i (-1)^i\beta^i(X,l)$.  Similarly,
 let $H_c^i(X,\Z_l)$ be the $l$-adic cohomology group
 of $X$ with compact support. Set $\beta^i_c(X,l):=\dim_{\Q_l}H_c^i(X,\Q_l)$
 the $i$-th $l$-adic Betti number of $X$ with compact support and
$\chi_c(X,l):=\sum_i (-1)^i\beta_c^i(X,l)$ the $l$-adic Euler-Poincar\'{e} Characteristic
with compact support. Note that $\chi_c(X,l)$ is independent of the
choice of $l$ prime to $char(K)$ (See, e.g., \cite{Katz} or \cite{Illusie}).

By using the method in the last section we will deduce  the  $l$-adic version of Lawson-Yau formula (see Theorem \ref{main}).

From the proof of equation (\ref{eq1.1}) above, we need  similar results for Lemma \ref{lemma2.1}-\ref{lemma2.3}, Proposition \ref{deform} and the homotopy invariance of
$l$-adic Euler-Poincar\'{e} Characteristics.

As we stated before, an algebraic version of Proposition \ref{deform} over any
algebraically closed field $K$ was proved  by Friedlander
(cf. \cite{Friedlander1}, Prop.3.2). Lemma \ref{lemma2.2} is a purely algebraic result
and hence it holds for any algebraically closed field $K$. An algebraic version of the
first statement in Lemma \ref{lemma2.3} follows from a corresponding result of
Lemma \ref{lemma2.1}. An algebraic version of the  second statement in
Lemma \ref{lemma2.3} holds over any algebraically closed field.
Therefore, the key part to prove Theorem \ref{main} is the following algebraic version of Lemma \ref{lemma2.1}
and the homotopy invariance of $l$-adic cohomology.

\begin{lemma} \label{lemma3.1}
Let $X$ be an algebraic variety (not necessarily smooth, compact, or irreducible) over an algebraically closed filed $K$ and let $Y\subset X$ be a closed algebraic set with complement $U$. Then $\chi(X, l)=\chi(U,l)+\chi(Y,l)$ for any positive integer $l$ prime to $char(K)$.
\end{lemma}
\begin{proof}
This lemma follows from the long localization exact sequence for $l$-adic cohomology and the following result proved by
Laumon (independently by Gabber).
\end{proof}

\begin{proposition}[\cite{Laumon}] For any algebraic variety $X$ over an algebraically closed field $K$ and integer $l$ prime to $char(K)$, we have
$\chi(X,l)=\chi_c(X,l)$.
\end{proposition}

To prove Theorem \ref{main}, we also need the following definition (cf. \cite{Friedlander1}, page 61).

\begin{definition}
A proper morphism $g:X'\to X$ of locally noetherian schemes is said to be a \textbf{bicontinuous algebraic morphism} if it is a set theoretic bijection and if for every $x\in X$ the associated map of residue fields $K(x)\to K(g^{-1}(x))$ is purely inseparable. A \textbf{continuous algebraic map} $f:X\to Y$ is  a pair ($g:X'\to X, f':X\to Y$) in which $g$ is a bicontinuous algebraic morphism (and $f$ is a morphism).
\end{definition}

\begin{lemma}
Let $F:X\times \A_K^1\to Y$ be a continuous algebraic map such that $F(-,0)=f(-)$ and $F(-,1)=g(-)$. Then the pullback $f^*:H^i(Y,\Z_l)\to H^i(X,\Z_l)$ is equal to $g^*: H^i(Y,\Z_l)\to H^i(X,\Z_l)$.
\end{lemma}
\begin{proof} It is essentially proved  in \cite{Friedlander1}, Prop. 2.1.
The map
$F:X\times \A_K^1\to Y$ induces a map in cohomology $F^*:H^i(Y,\Z_l)\to H^i(X\times \A_K^1,\Z_l)$. By the usual homotopy invariance of etale cohomology, one has the homotopy invariance  of $l$-adic cohomology, i.e., $H^i(X\times \A_K^1,\Z_l)\cong H^i(X,\Z_l)$. Therefore, we get the equality of $f$ and $g$ by the
restriction of $F$ to $X\times 0$ and $X\times 1$.
\end{proof}

\begin{corollary}\label{cor3.1}
Let $i:Y\subset X$ be an algebraically closed subset.
Let $F:X\times \A_K^1\to X$ be a continuous algebraic map such that $F(-,0)=id_X$, $F(x,1)\circ i=id_Y$ and
$F(y,t)=y$ for $y\in Y$. Then $i^*: H^i(X,\Z_l)\to H^i(Y,\Z_l)$ is an isomorphism.
\end{corollary}

The detailed computation is given below.

\begin{proof}[Proof of Theorem  \ref{main}]
Set $ \widetilde{Q}_{p,n}(t):=\sum_{d=0}^{\infty} \chi(C_{p,d}(\P^n))_K t^d,
$

$$T_{p+1,d}(\P^{n+1})_K:=\bigg\{ c=\sum n_iV_i\in C_{p+1,d}(\P^{n+1})_K| \dim(V_i\cap \P^n_K)=p, \forall i\bigg\}$$
and
$$B_{p+1,d}(\P^{n+1})_K:=\bigg\{ c=\sum n_iV_i\in C_{p+1,d}(\P^{n+1})_K| V_i\subset \P^n_K, \hbox{for some $i$}\bigg\}$$

By Corollary \ref{cor3.1} and the algebraic version of Proposition \ref{deform}, we have
\begin{equation}\label{eqn8}
 \chi(C_{p,d}(\P^n)_K)=\chi(T_{p+1,d}(\P^{n+1})_K).
\end{equation}

From the K\"{u}nneth formula for $l$-adic cohomology (cf. \cite{Milne}), the algebraic version of
Lemma \ref{lemma2.2} and Equation (\ref{eqn8}), we get

\begin{equation}\label{eqn9}
\chi(C_{p+1, d}(\P^{n+1})_K)=\chi(C_{p, d}(\P^{n})_K)+\sum_{i=1}^d \chi(C_{p+1, i}(\P^n)_K)\cdot\chi(C_{p,d-i}(\P^{n})_K),
\end{equation}
for integers $p\geq 0$ and $d\geq 1$.

By Corollary \ref{cor3.1}, we have $\chi(C_{0, d}(A^{n}_K))=1$, where $A^{n}_K$ is the
$n$-dimensional affine space over $K$. Hence the computation in the proof of Lemma \ref{lemma2.3}
works when $\C$ is replaced by any algebraically closed field $K$ and we
get the formula for $\chi(C_{0, d}(\P^{n})_K)$:
\begin{equation}\label{eqn10}
\chi(C_{0, d}(\P^{n})_K)=(^{n+1}_{d+1}).
\end{equation}

From the fact that $C_{p, d}(\P^{p+1})_K$ is
 the moduli space of hypersurfaces of degree $d$ in $\P^{p+1}_K$ and hence it
is a projective space over $K$ of dimension $(^{p+d+1}_{\quad d})-1$.
\begin{equation}\label{eqn11}
\chi(C_{p, d}(\P^{p+1})_K)=(^{p+d+1}_{\quad d}).
\end{equation}

From the definition of $\widetilde{Q}_{p,n}(t)$ and Equation (\ref{eqn9})-(\ref{eqn11}),
we get
\begin{equation}\label{eqn12}
\begin{array}{lll}
\widetilde{Q}_{p+1,n+1}(t)&=&\widetilde{Q}_{p+1,n}(t)\cdot \widetilde{Q}_{p,n}(t),\\
\widetilde{Q}_{0,m}(t)&=&(\frac{1}{1-t})^{m+1},\\
\widetilde{Q}_{q,q+1}(t)&=&(\frac{1}{1-t})^{q+2}.
\end{array}
\end{equation}

From Equation (\ref{eqn12}), we complete the proof of Theorem \ref{main} by induction on $n$.
\end{proof}

\section{Algebraic cycles with group action}
In this section, we apply our method to the space of algebraic cycles with certain finite group action $G$
to obtain the Euler characteristic of the $G$-invariant Chow varieties.
Let $G$ be a finite group of the automorphism of $\P^n$. By passing to an appropriate group
extension we can always assume that $\rho:G\to U_{n+1}$ and that its action on $\P^n$
 comes from the linear action of $U_{n+1}$ on homogeneous coordinates. The action of $G$ on $\P^n$ induces
 actions on the Chow varieties $C_{p,d}(\P^n)$.

 Denote by $$C_{p,d}(\P^n)^G:=\{c\in C_{p,d}(\P^n):g^*c=c, \forall g\in G\}$$
the $G$-invariant subset of $C_{p,d}(\P^n)$.  Since $G$ is a finite group of the automorphism of $\P^n$, it induces
an automorphism of $C_{p,d}(\P^n)$. Hence $C_{p,d}(\P^n)^G$ is a closed complex subvariety of $C_{p,d}(\P^n)$.

Choose homogeneous coordinates
$\C^{n+2}=\C^{n+1}\oplus\C$ for $\P^{n+1}=\Sigma\P^n$ and extend the fixed linear representation
$\rho:G\to U_{n+1}$  to a representation
 $\tilde{\rho}:G\to U_{n+2}$ by setting $\tilde{\rho}=\rho\oplus \lambda\cdot\id_{\C}$, where $\lambda\in \C^*$ is a fixed
 complex number. The construction was given in \cite{LM}, where $\lambda$ is chosen to be $1$.

Set $T_{p+1,d}(\P^{n+1})^G:=\big\{c=\sum n_iV_i\in C_{p+1,d}(\P^{n+1})^G|\dim(V_i\cap \P^{n})=p, ~ \forall i\big\}$
and $B_{p+1,d}(\P^{n+1})^G=C_{p+1,d}(\P^{n+1})^G-T_{p+1,d}(\P^{n+1})^G$.

The following proposition proved by Lawson and Michelsohn in \cite{LM} will be used in our calculation.
\begin{proposition}[\cite{LM}]
For each $p\geq 0$, $T_{p+1,d}(\P^{n+1})^G\subset C_{p+1,d}(\P^{n+1})^G$ is Zariski open.
 Moreover, the image $$\Sigma: C_{p,d}(\P^n)^G\to C_{p+1,d}(\P^{n+1})^G$$ is included in $T_{p+1,d}(\P^{n+1})^G$
 and $ \Sigma(C_{p,d}(\P^n)^G)\subset T_{p+1,d}(\P^{n+1})^G$ is a strong deformation retract.
\end{proposition}

In particular, their Euler characteristics coincide, i.e., we have
\begin{equation*}
 \chi(C_{p,d}(\P^n)^G)=\chi(T_{p+1,d}(\P^{n+1})^G).
\end{equation*}

\begin{remark}
In \cite{LM}, Lawson and Michelsohn consider the extension
$\tilde{\rho}=\rho\oplus \lambda\cdot\id_{\C}$ only for the case $\lambda=1$.
However, their proof  works for all  $\lambda\in\C^*$
without changing anything except that $1$ is replaced by $\lambda$ in $\tilde{\rho}$.
\end{remark}

From the definition,
$$B_{p+1,d}(\P^{n+1})^G=\big\{ c=\sum n_iV_i\in C_{p+1,d}(\P^{n+1})^G| V_i\subset \P^n, \hbox{for some $i$}\big\},$$ i.e.,
 there is at least one irreducible component lying in the  fixed $G$-invariant hyperplane $\P^n$.

\begin{lemma}\label{lemma4.02}
 $B_{p+1,d}(\P^{n+1})^G=\amalg_{i=1}^{d}B_{p+1,d}(\P^{n+1})_i^G$ , where $\amalg$ means disjoint union and
 $$B_{p+1,d}(\P^{n+1})_i^G=\left\{
                 c\in B_{p+1,d}(\P^{n+1})^G
		 \begin{array}{ll}\left|\begin{array}{ll}
		  &c=\sum n_kV_k+\sum m_jW_j,   \\
                  &  V_k\subset \P^n,\forall k,\\
		   & \dim(W_j\cap \P^n)=p, \forall j\\
		   & \deg(\sum n_kV_k)=i,\\
		   & and ~~\deg(\sum m_jW_j)=d-i.
		   \end{array}
                  \right.
                \end{array}
              \right\}.
$$
 For each $i$, $B_{p+1,d}(\P^{n+1})_i^G=C_{p+1, i}(\P^n)^G\times T_{p+1,d-i}(\P^{n+1})^G$.
\end{lemma}
\begin{proof}
An algebraic cycle $c\in B_{p+1,d}(\P^{n+1})_i^G$ may be written as
$c=\sum n_kV_k+\sum m_jW_j$ as the formal sum of irreducible varieties, where
$V_k\subset \P^{n}$ and $\dim(W_j\cap \P^{n})=p+1$. Since $c$ is $G$-invariant,
 we have  $\sum n_kV_k$ and $\sum m_jW_j$ are $G$-invariant. To see this, recall that
 $G$ is identified with the subgroup of  the unitary group $U_{n+1}$. Suppose  $g^*(V_i)=W_j$
 for some $g\in G$ and $i,j$. Since $\P^n$ is $G$-invariant, We have $\P^n=g^*(\P^n)\subset g^*(V_i)=W_j$.
This contradicts to the assumption that  $\dim(W_{j}\cap \P^{n})=p+1$. Similar for $\sum m_jW_j$.

Now the lemma follows  from the definition of $B_{p+1,d}(\P^{n+1})_i^G$.
\end{proof}

From Lemma \ref{lemma4.02}, we have $$\chi(B_{p+1,d}(\P^{n+1})_i^G)=\chi(C_{p+1, i}(\P^n)^G)\cdot\chi(T_{p+1,d-i}(\P^{n+1})^G).$$ Hence we get
$$
 \begin{array}{llll}
\chi(B_{p+1,d}(\P^{n+1}))
                 &=&\sum_{i=1}^d \chi(B_{p+1,d}(\P^{n+1})_i^G)\\
                 &&\hbox{(by inclusion-exclusion principle)}\\
                 &=&\sum_{i=1}^d \chi(C_{p+1, i}(\P^n)^G)\cdot\chi(T_{p+1,d-i}(\P^{n+1})^G)\\
                 &=& \sum_{i=1}^d \chi(C_{p+1, i}(\P^n)^G)\cdot\chi(C_{p,d-i}(\P^{n})^G). \\
                \end{array}
$$

Therefore we have the following result:

\begin{proposition}\label{prop4.03}
 For any integer $p\geq 0$ and $d\geq 1$, we have the following formula
\begin{equation}\label{eqn4.01}
\chi(C_{p+1, d}(\P^{n+1})^{G})=\chi(C_{p, d}(\P^{n})^G)+\sum_{i=1}^d \chi(C_{p+1, i}(\P^n)^G)\cdot\chi(C_{p,d-i}(\P^{n})^G),
\end{equation}
where $\chi(C_{q, 0}(\P^{N})^G)=1$ for integers $N\geq  q\geq 0$.
\end{proposition}

From our construction, we know that if the representation $\rho:G\subset U_{n+1}$ is
 \emph{diagonalizable}, i.e., up to a linear transformation,  $\rho=\oplus_{i=1}^{n+1} \lambda_i \cdot\id_{\C}$,
then equation (\ref{eqn4.01}) gives us a recursive formula. In these cases, the Euler
characteristic is calculated explicitly as follows.

\begin{theorem}\label{Th4.4}
Let  $\rho:G\subset U_{n+1}$ be a diagonalizable representation. The Euler characteristic of Chow variety of $G$-invariant
cycles $\chi(C_{p,d}(\P^{n})^G)$ is given by the formula
$$
\chi(C_{p,d}(\P^{n})^G)= \big(^{v_{p,n}+d-1}_{\quad\quad d}\big), \quad\hbox{where $v_{p,n}=(^{n+1}_{p+1})$}.
$$
\end{theorem}
\begin{proof}
The theorem follows from Proposition \ref{prop4.03} and the following initial values identities:
\begin{equation}
\chi(C_{0, d}(\P^{n})^G)=(^{n+d}_{~~d})
\end{equation}

As before, we can write $$C_{0,d}(\P^{n+1})^G=C_{0,d}(\C^{n+1})^G\coprod B_{0,d}(\P^{n+1})^G,$$
  where $C_{0,d}(\C^{n+1})^G\subset C_{0,d}(\P^{n+1})^G$ contains effective $G$-invariant 0-cycles $c$
  of degree $d$ such that no points in $c$ lying in the fixed hyperplane $\P^n$ and
   $B_{0,d}(\P^{n+1})^G$ is the complement of $C_{0,d}(\C^{n+1})^G$ in  $C_{0,d}(\P^{n+1})^G$.
   We claim that $C_{0,d}(\C^{n+1})^G$ is contractible. To see this, note that $\P^{n+1}-\P^{n}=\C^{n+1}$ and
   Let $\phi_t:\C^{n+1}\to \C^{n+1}$ denote scalar multiplication by $t\in\C$.
   The family of maps $\phi_t$ induces a family of maps $\phi_{t*}:C_{0,d}(\C^{n+1})\to C_{0,d}(\C^{n+1})$
   since the multiplication by $t\in\C$ is $G$-invariant. From the definition, the map $\phi_{1*}=\id$ and $\phi_{0*}$ is a constant map.

   We can write $$B_{0,d}(\P^{n+1})^G=\coprod_{i=1}^d B_{0,d}(\P^{n+1})_i^G$$ as in Lemma \ref{lemma4.02}, where
 $B_{0,d}(\P^{n+1})_i^G$ contains $G$-invariant 0-cycles $c$ of degree $d$
 on $\P^{n+1}$ in which there are exact $i$ points (count multiplicities) lying in $\P^n$, hence
 $B_{0,d}(\P^{n+1})_i^G=C_{0,i}(\P^n)^G\times C_{0,d-i}(\C^{n+1})^G$. In particular,
 $\chi(B_{0,d}(\P^{n+1})^G)=\sum_{i=1}^d\chi(C_{0,i}(\P^n)^G)$.

 Therefore, we have
 $$\chi(C_{0, d}(\P^{n+1})^G)=1+\sum_{i=1}^d\chi(C_{0,i}(\P^n)^G).$$
Now the  formula in the lemma follows from this by induction.

\end{proof}

\begin{remark}
By a carefully checking the proof of Theorem 4.6 in \cite{LM} and the proof of Theorem \ref{Th4.4} above, we observe that if the linear representation
$\rho:G\to U_{n+1}$ is diagonalizable, then the conclusion in Theorem \ref{Th4.4} holds
even if there is no assumption of finiteness of $G$.
More precisely, let $T_{n+1}\subset U_{n+1}$ be the maximal torus and let $G$ be any subgroup of $T_{n+1}$. Then
 The Euler characteristic of Chow variety of $G$-invariant
cycles $\chi(C_{p,d}(\P^{n})^G)$ is given by the formula
$$
\chi(C_{p,d}(\P^{n})^G)= \big(^{v_{p,n}+d-1}_{\quad\quad d}\big), \quad\hbox{where $v_{p,n}=(^{n+1}_{p+1})$}.
$$
This explains Theorem 4.1, one of the main results in \cite{Lawson-Yau}, on the invariance of Euler characteristic of
a compact complex analytic space and that of the fixed-point set under a holomorphic $S^1$-action in the  important case
for Chow varieties over $\C$.

\end{remark}

\section{The Euler Characteristic for the space of right-quaternionic cycles }
Let $\H$ denote the quaternions with standard basis $1,\textbf{i},\textbf{j},\textbf{k}$,
and let $\C^2\stackrel{\cong}{\to} \H$ be the canonical isomorphism
given by $(u,v)\mapsto u+v\textbf{j}$.
This gives us a canonical complex isomorphism $\C^{2n}\stackrel{\cong}{\to} \H^n$.

Under this identification \emph{right} scalar multiplication
by $\textbf{j}$ in $\H^n$ becomes the complex linear map
$$J:\C^{2n}\to \C^{2n}, \quad J(u_1,...,u_n,v_1,...,v_n)=(-v_1,...,-v_n,u_1,...,u_n).$$

This induces a holomorphic map $\bar{J}:\P^{2n-1}\to\P^{2n-1}$ with
$\bar{J}=Id$. Note that the fixed point set of $\bar{J}$ is a pair of
 disjoint $\P^{n-1}$. The involution $\bar{J}$ carries algebraic subvarieties
 of $\P^{2n-1}$ to themselves and induces a \emph{holomorphic}  involution
\begin{equation}\label{eqn4.0}
\bar{J}_*:C_{p,d}(\P^{2n-1})\to C_{p,d}(\P^{2n-1})
\end{equation}
 for all $p$ and $d$.

\begin{remark}
The construction is an analog to the one given in \cite{LLM}, where  a \emph{left} scalar
multiplication by $j$ on $\H^n$ was checked in detail. We consider the \emph{right} multiplication
here since the induced map $\bar{J}_*$ is a holomorphic.
\end{remark}

Let $C_{p,d}(n)\subset C_{p,d}(\P^{2n-1})$ denote the $\bar{J}_*$-fixed point set, i.e.,
the set of $\bar{J}_*$-invariant algebraic $p$-cycles. An element $c\in C_{p,d}(n)$ is called a \textbf{right
quaternionic cycle}. Since $\bar{J}_*$
is a holomorphic involution, $C_{p,d}(n)$ is a closed complex algebraic set.

Since $J$ is diagonalizable, in fact, $$J\sim diag \left\{
\begin{array}{ccc}
\left(\begin{array}{ccc}
\sqrt{-1}&0\\
0&-\sqrt{-1},
\end{array}\right),
&\cdots, &
\left(\begin{array}{ccc}
\sqrt{-1}&0\\
0&-\sqrt{-1}
\end{array}\right)
\end{array}\right\}.
$$

As an application of Theorem \ref{Th4.4}, we have the following result.
\begin{corollary}\label{cor5.2}
 For any $p\geq 0$, we have
$$\chi(C_{p,d}(n))=\chi(C_{p,d}(\P^{2n-1}))=\big(^{v_{p,2n-1}+d-1}_{\quad\quad d}\big), \quad\hbox{where $v_{p,2n-1}=(^{2n}_{p+1})$}.$$
\end{corollary}

\begin{example}
 For $p=0$, we have $\chi(C_{p,d}(n))=(^{2n+d-1}_{\quad d})$.
 \end{example}
Alternatively, this can be seen in the following way.
The set $C_{0, d}({n})$ can be decomposed into the disjoint union of quasi-projective algebraic varieties according to the number
of fixed points of $\bar{J}$ lying in $\P^{n-1}\amalg \P^{n-1}$, i.e.,
 $C_{0, d}({n})=\amalg_{i=0}^{d} \sp^i(\P^{n-1}\amalg \P^{n-1})\times \sp^{\frac{1}{2}(d-i)}(G(n))$,
where $G(n)=\P^{2n-1}-(\P^{n-1}\amalg \P^{n-1})$ is a bundle over $\P^{n-1}$ with fibers $\C^{n}-\{0\}$.
Hence we have $$\chi(C_{0, d}({n}))=\sum_{i=0}^{d}\chi(\sp^i(\P^{n-1}\amalg \P^{n-1})\cdot\chi(\sp^{\frac{1}{2}(d-i)}(G(n)).$$
Since the Euler characteristic of $G(n)$ is zero, we get $\chi(\sp^m(G(n)))=0$ for all $m>0$ by MacDonald formula (cf. \cite{Macdonald}).
Therefore, $\chi(C_{0, d}({n}))=\chi(\sp^d(\P^{n-1}\amalg \P^{n-1})).$  Note that
\begin{equation}\label{eqn4.04}
\sp^d(\P^{n-1}\amalg \P^{n-1})=\amalg_{i=0}^{d}\sp^i(\P^{n-1})\times \sp^{d-i}(\P^{n-1}),
\end{equation}
where a 0-cycle $c=\sum n_i P_i+\sum m_j P_j'\in \sp^d(\P^{n-1}\amalg \P^{n-1})$ of degree $d$ is written as the sum of two 0-cycles such that
$P_i$ is in the first copy of $\P^{n-1}$ but $P_j'$ is in the second copy of $\P^{n-1}$.
By equation (\ref{eqn4.04}) and the fact that $\chi(C_{0,i}(\P^{n-1}))=(^{n+i-1}_{\quad i})$, we get
$$
\begin{array}{lll}
\chi(\sp^d(\P^{n-1}\amalg \P^{n-1}))
&=&\sum_{i=0}^{d}\chi(\sp^i(\P^{n-1}))\cdot \chi(\sp^{d-i}(\P^{n-1}))\\
&=&\sum_{i=0}^{d}(^{n+i-1}_{\quad i})\cdot (^{n+d-i-1}_{\quad d-i})\\
&=& (^{2n+d-1}_{\quad d}),
\end{array}
$$
where the last equality is obtained by comparing the coefficients of $t^d$ in the Taylor series of
$\frac{1}{(1-t)^{2n}}=\frac{1}{(1-t)^{n}}\cdot\frac{1}{(1-t)^{n}}$.

\begin{example}
For $d=1$, we have $\chi(C_{p,1}(n))=(^{2n}_{p+1})$.
\end{example}
Alternatively, this can be seen in the following way.
The eigenvalues of $J$ are $\pm \sqrt{-1}$, each of them is of multiplicity $n$.
Let $\{e_i\}_{1\leq i\leq n}$ be the eigenvectors of the eigenvalue $\sqrt{-1}$ and
$\{f_i\}_{1\leq i\leq n}$ be the eigenvectors of the eigenvalue $-\sqrt{-1}$.
The $J$-invariant $(p+1)$-complex vector space is spanned by $i$ eigenvectors from
$\{e_i\}_{1\leq i\leq n}$ and $p+1-i$ eigenvectors from $\{f_i\}_{1\leq i\leq n}$. Therefore,
$$
C_{p,1}(n)=C_{p,1}(\P^{2n-1})^{\bar{J}_*}=\coprod_{1\leq i\leq p+1} G(i, n)\times G(p+1-i,n),
$$
where $G(i,n):=G(i,\C^{n})$ is the Grassmannian of $i$-dimensional complex linear subspaces in $\C^{n}$.
Therefore,
$$
\begin{array}{lll}
\chi(C_{p,1}(n)) =\chi(C_{p,1}(\P^{2n-1})^{\bar{J}_*})
&=&\sum_{i=1}^{p+1}  \chi(G(i, n))\cdot \chi(G(p+1-i,n)),\\
&=&  \sum_{i=1}^{p+1}   (^{n}_{ i})\cdot (^{\quad n}_{ p+1-i})\\
&=&  (^{2n}_{p+1}),
\end{array}
$$
where the last equality is obtained by comparing the coefficients of $t^{p+1}$ in the binomial expansion of
${(1+t)^{2n}}={(1+t)^{n}}\cdot{(1+t)^{n}}$.

\end{document}